\documentclass[12pt]{article}
\usepackage{amssymb,amsmath}

\def\={~=~}
\newcommand{\lp}{\left(}

\newcommand{\rp}{\right)}

\newcommand{\N}{\mathbb{N}}

\newcommand{\Z}{\mathbb{Z}}

\begin{document}
\thispagestyle{empty}

\begin{center}
\Large
Trigonometric Series via Laplace Transforms
\end{center}

\vspace*{0.2 in}

\begin{flushright}
Costas J. Efthimiou\footnote{costas@physics.ucf.edu}  \\
Department of Physics \\
University of Central Florida \\
Orlando, FL 32816  \\
USA
\end{flushright}

\paragraph*{Introduction}
In another NOTE in this MAGAZINE, the author presented a method
\cite{Efthimiou} that uses the Laplace transform and allows one to
find exact values for a large class of convergent series of
rational terms. Recently, in this MAGAZINE too, Lesko and Smith
\cite{Lesko} revisited the method and demonstrated an extension of
the original idea to additional infinite series. Our intention in
this note is to illustrate the power of the technique in the case
of trigonometric series.

\paragraph*{Trigonometric series}
Trigonometric series play a vital role in mathematics and physics.
Many results are known but most of them can be obtained only
via Fourier analysis. For example
\begin{eqnarray}
  \sum_{n=1}^\infty {\cos(nx)\over n} &=& -\ln\lp2\sin{x\over2}\rp~,
  ~~~~~0<x<2\pi~, \label{eq:1}\\
  \sum_{n=1}^\infty {\cos(nx)\over n^2} &=& {3x^2-6\pi x +2\pi^2\over 12}~,
   ~~~~~0\le x \le 2\pi~, \label{eq:2}\\
  \sum_{n=1}^\infty {\sin(nx)\over n} &=& {\pi-x\over2}~,
  ~~~~~0<x<2\pi~, \label{eq:3}\\
  \sum_{n=1}^\infty (-1)^n{\cos(nx)\over n} &=& -\ln\lp2\cos{x\over2}\rp~,
  ~~~~~-\pi<x<\pi~,\\
   \sum_{n=1}^\infty (-1)^n{\cos(nx)\over n^2} &=& {3x^2-\pi^2\over 12}~,
   ~~~~~-\pi\le x \le \pi~,\\
  \sum_{n=0}^\infty {\cos((2n+1)x)\over (2n+1)} &=& -{1\over2}
       \ln\tan{x\over2}~,
  ~~~~~0<x<\pi~,\\
  \sum_{n=0}^\infty {\cos((2n+1)x)\over (2n+1)^2} &=&
       {\pi^2-2\pi x\over8}~,
  ~~~~~0\le x\le \pi~.
\end{eqnarray}
These and many more results may be found in \cite{Davies}. Given a
function, it is relatively straightforward to expand it in a
trigonometric series. However, it is almost impossible to guess
the function that would generate a given trigonometric series when
as its Fourier series. For example, given the series
$$
   \sum_{n=0}^\infty (-1)^n \,{\sin((2n+1)x)\over (2n+1)^2}~,
   ~~~-\pi/2\le x\le \pi/2~,
$$
it is not easy to guess that the function
$$
  f(x) \= {\pi x\over4}
$$
will give the result sought and then proceed to prove it. On the
other hand, the method of \cite{Efthimiou},\cite{Lesko} can obtain
the results in a straightforward manner with no ad hoc guessing.

\paragraph*{The Method}
As Lesko and Smith have pointed out, the original method of \cite{Efthimiou}
can be applied to series of the form $\sum_{n\in I} u_n v_n$ where
$I$ is a subset of $\Z$. In series of this form, it is often convenient
to write only one of the factors, say $v_n$, as a Laplace transform
of a function $f(t)$
$$
    v_n \= \int_0^{+\infty} e^{-nt} \, f(t)~.
$$
Then
\begin{eqnarray*}
  \sum_{n\in I} u_n v_n
         &=& \sum_{n\in I}u_n\,\int_0^{+\infty}e^{-nt}\,f(t)\,dt~.
\end{eqnarray*}
Assuming that the order of the operations of summation and integration
can be exchanged
\begin{eqnarray*}
  \sum_{n\in I} u_n v_n
         &=& \int_0^{+\infty}\lp\sum_{n\in I}u_n\,e^{-nt}\rp\,f(t)\,dt~.
\end{eqnarray*}
In this note we shall always exchange the order of the two operations
assuming that the reader knows to reason for its validity. One may
consult \cite{Efthimiou},\cite{Lesko} for details.
If one can find an explicit
function $h(t)=\sum_{n\in I}u_n\,e^{-nt}$,
then she has succeeded to write the initial series in a
simple integral representation:
\begin{eqnarray*}
    \sum_{n\in I} u_n v_n &=& \int_0^{+\infty} h(t)\,f(t)\,dt~.
\end{eqnarray*}
If, furthermore, the integration
can be performed, then analytic answers for the initial series are
obtained.

\paragraph*{Easy trigonometric sums}
To apply the method to trigonometric series, we need to be able to handle
series of the form
\begin{eqnarray*}
   S &=& \sum_{n\in I} \sin(nx)\,e^{-nt}~,\\
   C &=& \sum_{n\in I} \cos(nx)\,e^{-nt}~.
\end{eqnarray*}
These summations are performed easily using complex-number notation:
$$
  C+iS \=  \sum_{n\in I} e^{inx}\, e^{-nt}~.
$$
In particular, when $I=\N^*=\{1,2,3,...\}$, assuming that $x$ is a
real number and $t>0$,
\begin{eqnarray*}
   \sum_{n=1}^\infty \sin(nx)\,e^{-nt}
        &=& {e^{-t}\sin x\over 1-2\cos x\, e^{-t}+e^{-2t}}~,\\
   \sum_{n=1}^\infty \cos(nx)\,e^{-nt}
        &=& {e^{-t}(\cos x-e^{-t})\over 1-2\cos x\, e^{-t}+e^{-2t}}~.
\end{eqnarray*}
Similarly we can compute other sums. Also, it is possible to start
with these sums, and by changes in the argument $x$ and simple
manipulations, derive formul\ae\ for new sums, such as sums over
the even or odd integers only. The reader may wish to experiment
with this idea.

\paragraph*{Trigonometric series via the Laplace transform}
We are now ready to find exact sums for more complicated
trigonometric series. We shall demonstrate the method with two
examples, namely (\ref{eq:1}) and (\ref{eq:3}). The other
formul\ae\ may be obtained similarly. Alternatively, they may be
computed using algebraic and integral operations on (\ref{eq:1})
and (\ref{eq:3}). For example, (\ref{eq:2}) may be proved by
integrating (\ref{eq:3}) and using the well known sum
$$
   \sum_{n=1}^\infty {1\over n^2} ~=~ {\pi^2\over 6}~.
$$
The last sum, in turn, may be obtained easily using the original
method \cite{Efthimiou} or various other techniques.

The approach described here, however, allows us to derive equation
(\ref{eq:2}) without reference to any other sum, assuming that one
can integrate the necessary functions, perhaps using a table of
integrals. Since tables of integrals are widely available and they
are quite extensive (for instance, \cite{ref4}), the method seems
to be quite effective and straightforward.

\begin{itemize}
\item
We start with the series
$$
  \sum_{n=1}^\infty {\cos(nx)\over n^\nu}
$$
where $\nu\in\N^*$ and $0<x<2\pi$ if $\nu=1$ or $0\le x\le 2\pi$
if $\nu>1$. Using
$$
  {1\over n^\nu} \= {1\over(\nu-1)!}\,
   \int_0^\infty e^{-nt} \, t^{\nu-1} \, dt~,
$$
we write
\begin{eqnarray*}
 \sum_{n=1}^\infty {\cos(nx)\over n^\nu}
  &=& {1\over(\nu-1)!}\,
      \sum_{n=1}^\infty \cos(nx) \int_0^\infty e^{-nt}\, t^{\nu-1}\, dt \\
  &=& {1\over(\nu-1)!}\,
  \int_0^\infty \lp\sum_{n=1}^\infty \cos(nx) e^{-nt}\rp\, t^{\nu-1}\, dt \\
  &=& {1\over(\nu-1)!}\,
    \int_0^\infty  {e^{-t}(\cos x-e^{-t})\over 1-2\cos x\, e^{-t}+e^{-2t}}
      \, t^{\nu-1}\, dt~.
\end{eqnarray*}
With the change of variables $u=e^{-t}$ the integral is cast in a more
compact form
$$
  \sum_{n=1}^\infty {\cos(nx)\over n^\nu}
  \= {(-1)^{\nu-1}\over(\nu-1)!}\,
      \int_0^1  {\cos x-u\over 1-2\cos x\, u+u^2}
      \, (\ln u)^{\nu-1}\, du~.
$$
For $\nu=1$
\begin{eqnarray*}
 \sum_{n=1}^\infty {\cos(nx)\over n}
  &=& -{1\over2}\, \int_0^1  {d(1-2\cos x\,u+u^2)\over 1-2\cos x\, u+u^2}\\
  &=& -{1\over2}\, \ln(1-2\cos x\, u+u^2)\Big|_0^1
  \= -\ln(2\sin{x\over2})~.
\end{eqnarray*}

\item
Following the same steps for the series
$$
  \sum_{n=1}^\infty  {\sin(nx)\over n^\nu}~,
$$
we will arrive at the integral representation
$$
  \sum_{n=1}^\infty {\sin(nx)\over n^\nu}
  \= {(-1)^{\nu-1}\over(\nu-1)!}\,\sin x \,
      \int_0^1  {(\ln u)^{\nu-1}\over 1-2\cos x\, u+u^2}
      \, du~.
$$
In particular, for $\nu=1$
\begin{eqnarray*}
 \sum_{n=1}^\infty {\sin(nx)\over n}
  &=& \sin x\, \int_0^1  {1\over (u-\cos x)^2+\sin^2x}\,du\\
  &=& \tan^{-1}{u-\cos x\over\sin x}\Big|_0^1
  \= \tan^{-1}{\sin x\over1-\cos x}\={\pi-x\over2}~.
\end{eqnarray*}

\end{itemize}

\paragraph*{Conclusion}
Although the results presented in this paper are not new, we hope
that the reader will appreciate the ease and transparency of the
method. Given \textit{any} series such as those described in the
original articles of \cite{Efthimiou}, \cite{Lesko}, and this
note, the steps are well-defined and require no special tricks.
However, traditional methods do vary from series to series, tricks
may be necessary to be introduced, and some (or a lot) guessing
might be required. We invite the reader to verify our claim by
using the Laplace transform technique to find exact sums for her
favorite series (of the type described in \cite{Efthimiou},
\cite{Lesko}, and this note) and then compare with the traditional
methods.

\paragraph*{Note Added in Proof}
While this article was under review, we received a message from
Harvey J. Hindin, who pointed out that Albert D. Wheelon had also
used the idea of integral transformations to compute exact values
for infinite series \cite{ref5}.


\end{document}